\documentclass[a4, 12pt]{amsart}
\usepackage[mathscr]{eucal}
\usepackage{amssymb}
\usepackage{latexsym}
\usepackage{amsthm}
\usepackage{color}
\theoremstyle{plain}
\newtheorem{theorem}{Theorem}[section]

\newtheorem{corollary}{Corollary}[section]

\newtheorem{remark}{Remark}[section]
\newtheorem{lemma}{Lemma}[section]
\newtheorem{proposition}{Proposition}[section]

\makeatletter
\@addtoreset{equation}{section}

\title[Translators and expanders with constant curvature]
{Rigidity and classification of complete translating solitons and self-expanders with constant curvature in $\mathbb{R}^{4}$}

\author [Z. Li and C. Xing]{Zhi Li and Cheng Xing}

\address{Zhi Li \\  \newline \indent School of Mathematics and Statistics, Henan Normal University,
\newline \indent 453007, Xinxiang, Henan, China.}
\email{lizhihnsd@126.com}

\address{Cheng Xing\\  \newline \indent School of Mathematics and Statistics, Henan Normal University,
\newline \indent 453007, Xinxiang, Henan, China.}
\email{xingcheng@htu.edu.cn}

\begin{document}
\maketitle

\begin{abstract}
For complete immersed surfaces in Euclidean four-space whose mean curvature
vector has constant length, translators are affine planes parallel to their fixed
nonzero velocity, and self-expanders are planes through their centre. This rigidity
requires no Lagrangian assumption, properness, volume growth control, or bound
on the second fundamental form. Complete Lagrangian translators and self-expanders in $\mathbb{C}^{2}$ with constant Gaussian curvature are flat and, up to unitary
equivalence, are properly embedded products of two complete planar curves.
Translating factors are lines or grim reapers. Expanding factors are lines through
the origin or nonstraight self-expanders. The proofs use moving frames, a
transnormal function, a cubic tensor identity, and analytic continuation.
\end{abstract}

\footnotetext{2020 \textit{Mathematics Subject Classification}:Primary 53E10; Secondary 53C40, 53D12.}
\footnotetext{{\it Key words and phrases}: Mean curvature flow, Translating soliton, Self-expander, Lagrangian surface, Gaussian curvature.}

\section{introduction}
\vskip2mm
\noindent

Mean curvature flow evolves a family of immersed submanifolds $x_{t}=x(\cdot, t): M^{n}\to \mathbb{R}^{n+p}$ by
\begin{equation*}
\frac{d}{dt}x(p, t)=\vec{H}(p,t), \ \ x(p, 0)=x(p), \ \ p\in M^{n},
\end{equation*}
where $\vec{H}(p,t)$ is the mean curvature vector of $M_{t}$ at $p\in M^{n}$. Its self-similar solutions play a central role both in the local analysis of singularities and in the long-time behaviour of the flow. Two basic examples are translating solitons (or, simply, translators), characterised by
\begin{equation}\label{1.1}
\vec{H}=V^{\perp},\quad 0\neq V\in\mathbb{R}^{n+p}\text{ constant},
\end{equation}
and self-expanders, satisfying
\begin{equation}\label{1.2}
\vec{H}=x^{\perp}.
\end{equation}
Under our convention $\Delta x=\vec{H}$, equation \eqref{1.1} generates translation by $tV$, while \eqref{1.2} generates the homothetic flow $x_{t}=\sqrt{2t}x$ for $t>0$, up to tangential reparametrisation. The normalisation in \eqref{1.2} is essential for the numerical coefficients below. Basic examples and the theory of translators and expanders may be found in \cite{AltschulerWu1994,AncariCheng2022,ChengZhou2018,Clutterbuck2007,Halldorsson2012,Halldorsson2013,Smoczyk2021}.

The interaction between intrinsic and extrinsic geometry becomes particularly interesting for Lagrangian mean curvature flow in $\mathbb{C}^{n}$. The Lagrangian condition is preserved along smooth mean curvature flow, and the parallel complex structure identifies the tangent and normal bundles. This produces additional symmetry in the second fundamental form, a feature absent for general submanifolds. Foundational constructions and rigidity results include \cite{Anciaux2006,CastroLerma2010,CastroLerma2014,ChengHoriWei2022,JoyceLeeTsui2010,LeeWang2009,NevesTian2013}.

From view points of submanifolds theory, it is very natural to study rigidity and classification theorems for the Lagrangian self-shrinkers, the Lagrangian
translators and the Lagrangian self-expanders. Two curvature assumptions that often arise in classification problems must be distinguished. The extrinsic condition constant $|\vec{H}|^{2}$ fixes the mean curvature length. The intrinsic condition constant $K$ fixes the Gaussian curvature. For an immersed surface in $\mathbb{R}^{4}$, the Gauss equation gives
$2K=|\vec{H}|^{2}-|B|^{2}$. Thus neither of these two constancy assumptions implies the other without further information. Recent results on Lagrangian
translators and Lagrangian self-expanders under hypotheses involving $|\vec{H}|^2$ or $|B|^{2}$ appear in \cite{LiLiuQiao2021,LiWei2024,LiWei2025}; complete Lagrangian self-shrinkers are investigated in \cite{ChengHoriWei2022,LiWangWei2025,LiWeiXing2026}. In contrast, our second classification problem below concerns constant Gaussian curvature.

We first address the mean curvature length problem without imposing the Lagrangian condition.

\begin{theorem}\label{theorem 1.1}
Let $x:M^{2}\to\mathbb{R}^{4}$ be a smooth, complete, connected translator without boundary. If $|\vec{H}|^2$ is constant, then $\vec{H}=0$ and $B=0$.
Consequently, $x$ is a global isometry onto an affine two-plane $p_{0}+P$ with $V\in P$.
\end{theorem}

\begin{theorem}\label{theorem 1.2}
Let $x:M^{2}\to\mathbb{R}^{4}$ be a smooth, complete, connected self-expander without boundary. If $|\vec{H}|^2$ is constant, then $\vec{H}=B=0$, and $x$ is a global isometry onto a linear two-plane through the origin.
\end{theorem}

\begin{remark}
The above conclusion, without conditions the squared norm of the second fundamental form being bounded from above and Lagrangian condition, proves that the conclusion of \cite{LiWei2025} still holds true. The same transnormal function occurs for self-shrinkers, with a different
Laplacian, in the preprint of Li, Wei and Xing \cite{LiWeiXing2026}.
\end{remark}
\begin{remark}
The mechanisms are different. For translators, the nonzero mean curvature case is ruled out locally using a moving frame, the Gauss equation, and Codazzi compatibility. For expanders, completeness and a transnormal identity force the constant mean curvature to vanish. In particular, the translator's vanishing conclusion does not itself require completeness; completeness is used only to identify the entire plane.
\end{remark}

Our second condition prescribes the intrinsic Gaussian curvature $K$.
In particular, constant $K$ fixes the difference of two extrinsic curvature functions. The
constant $|B|^{2}$ problem for complete Lagrangian surfaces was treated by Li, Liu and Qiao for translators
\cite{LiLiuQiao2021} and by Li and Wei for expanders \cite{LiWei2024}. Chen \cite{Chen2005} classified general
constant-curvature Lagrangian surfaces in $\mathbb{C}^{2}$ locally on an open dense set.
For the intrinsic problem, let $J$ be the standard complex structure on $\mathbb{C}^2$. An immersion is Lagrangian when $J(TM)=NM$. The stronger cubic-tensor symmetry leads to global classifications.

\begin{theorem}\label{theorem 1.3}
Let $x:M^{2}\to\mathbb{C}^2$ be a smooth, complete, connected Lagrangian translator without boundary. If If the Gaussian curvature is constant, then the surface is flat and isometric to $\mathbb{R}^{2}$, and there exist a unitary transformation $U$ and complete unit-speed plane curves $\gamma_{i}:\mathbb{R}\to\mathbb{C}$ such that
\begin{equation*}
x(s,t)=U\big(\gamma_{1}(s),\gamma_{2}(t)\big),\quad
\gamma^{''}_{i}=(V_{i})^{\perp_{\gamma_{i}}},\quad (V_{1},V_{2})=U^{-1}V.
\end{equation*}
Here $\perp_{\gamma_{i}}$ denotes projection onto the normal line of the curve. Each factor is a straight line or a grim reaper curve; the product is a proper embedding. Conversely, every such product is a complete flat Lagrangian translator, provided the factor equations hold for the prescribed $V$.
\end{theorem}

\begin{theorem}\label{theorem 1.4}
Let $x:M^{2}\to\mathbb{C}^2$ be a smooth, complete, connected Lagrangian self-expander without boundary. If its gaussian curvature $K$ is constant,, then $K=0$, and $M$ is isometric to $\mathbb{R}^{2}$, and
\begin{equation*}
x(s,t)=U\bigl(\gamma_1(s),\gamma_2(t)\bigr),\quad
\gamma^{''}_{i}=(\gamma_{i})^{\perp_{\gamma_{i}}},
\end{equation*}
for a unitary transformation $U$ and complete unit speed planar self-expanders $\gamma_{i}$. Each factor is a line through the origin or a nonstraight complete expanding curve. The product is properly embedded. Conversely, each such product is a complete flat Lagrangian self-expander.
\end{theorem}

\begin{remark}
In fact, the authors of
the present paper \cite{LiWeiXing202611} proved that the conclusion remains valid for self-shrinkers under the same hypotheses as those in Theorems \ref{theorem 1.3} and \ref{theorem 1.4}.
\end{remark}

Joyce, Lee and Tsui \cite{JoyceLeeTsui2010} constructed families of Lagrangian translators and self-expanders. The translating products of a line and a grim reaper, or of two grim reapers,
already appear in Castro and Lerma \cite{CastroLerma2012}. For self-shrinkers, Li and Wang \cite{LiWang2017} obtained
a product decomposition for complete flat Lagrangian surfaces. Here the
exclusion of negative constant $K$ uses an identity for the covariant derivative of the Lagrangian cubic
tensor. After proving flatness, we extend the local product decomposition to the universal cover. The
known classification of planar soliton curves \cite{Halldorsson2012} then gives injectivity and properness on
the original surface.

The exclusion of negative constant Gaussian curvature rests on an explicit identity for the covariant derivative of the Lagrangian cubic tensor, combined with an approximate extremal sequence. Flatness is then promoted to a global Cartesian splitting by commutation of shape operators and real-analytic elliptic regularity. Finally, planar curve classification gives embeddedness and eliminates nontrivial deck transformations. We include the calculations and the completeness arguments to make precise where each hypothesis enters.

The paper is organised as follows. Section \ref{section 2} fixes conventions and gives auxiliary lemmas. Sections~\ref{section 3}-\ref{section 4} prove the general-surface rigidity statements. Sections~\ref{section 5}-\ref{section 7} establish the curvature identities and exclude nonzero constant $K$. Sections~\ref{section 8}-\ref{section 9} prove the global product classification.

\vskip5mm
\section {Preliminaries}\label{section 2}
\vskip2mm

\noindent

Let $x:M^{2}\to\mathbb{R}^{4}$ be a smooth immersion. Throughout, $M$ is connected and without boundary, and completeness refers to the induced metric. For an orthonormal tangent
frame  $\{e_{1},e_{2}\}$, write $B_{ij}=B(e_{i},e_{j})$. Thus
\begin{equation}\label{2.1}
\vec{H}=\operatorname{tr}_{g} B=\sum_{i=1}^{2} B_{ii},\quad \Delta x=\vec{H},
\end{equation}
where $B$ is the second fundamental form. Denote by $D$, $\nabla$ and $\nabla^{\perp}$ the Euclidean, Levi-Civita, and normal connections, respectively. For tangent vector fields $X,Y$ and a normal field $\eta$, the Gauss and Weingarten formulas are
\begin{equation}\label{2.2}
D_{X}Y=\nabla_{X}Y+B(X,Y), \quad D_{X}\eta=-A_{\eta}X+\nabla^{\perp}_{X}\eta,
\end{equation}
where $\langle A_{\eta}X, Y\rangle=\langle B(X,Y), \eta\rangle$. Our curvature convention is $$R(X,Y)=\nabla_{X}\nabla_{Y}-\nabla_{Y}\nabla_{X}-\nabla_{[X,Y]}.$$
With $K=\langle R(e_{1},e_{2})e_{2}, e_{1}\rangle$, the Gauss equation gives
\begin{equation}\label{2.3}
K=\langle B_{11},B_{22}\rangle-|B_{12}|^{2},
\quad 2K=|\vec{H}|^{2}-|B|^{2},
\end{equation}
The covariant derivative of $B$ is
$$(\nabla_{X}B)(Y,Z)=\nabla_X^{\perp}(B(Y,Z))-B(\nabla_{X}Y,Z)-B(Y,\nabla_{X}Z),$$
and the Euclidean Codazzi equation is
\begin{equation}\label{2.4}
(\nabla_{X}B)(Y,Z)=(\nabla_{Y}B)(X,Z).
\end{equation}

\begin{lemma}\label{lemma 2.1}(\cite{LiWeiXing2026})
Let $(N^{m},g)$ be a nonempty complete Riemannian manifold. If $u$ is smooth and nonnegative and satisfies $|\nabla u|^2=4u$. Then $u$ vanishes somewhere.
\end{lemma}
\begin{proof}
If $u\equiv0$, the conclusion is immediate. On $\Omega=\{u>0\}$ define $r=\sqrt{u}$, then $|\nabla r|=1$. Fix $p\in\Omega$ and take the maximal integral curve $\gamma:[0,\tau)\to\Omega$ of $-\nabla r$ with $\gamma(0)=p$. Its speed equals one and
$$
\frac{d}{dt}r(\gamma(t))=-1,\quad r(\gamma(t))=r(p)-t,
$$
so $\tau\le r(p)$. The curve is Cauchy as $t\uparrow\tau$ because $d(\gamma(t),\gamma(s))\leq|t-s|$. Completeness supplies $q=\lim_{t\uparrow\tau}\gamma(t)$. If $\tau<r(p)$, then $u(q)>0$ and the smooth vector field $-\nabla r$ extends the integral curve, contradicting maximality. Hence $\tau=r(p)$ and $u(q)=0$.
\end{proof}

\begin{lemma}\label{lemma 2.2}
Under the hypotheses of Lemma \ref{lemma 2.1}, if $u(q)=0$, the
self-adjoint Hessian endomorphism  $\mathcal{A}$ defined by $\langle\mathcal{A}X, Y\rangle=\nabla^{2}u(X,Y)|_{q}$
satisfies $\mathcal{A}^{2}=2\mathcal{A}$. Consequently, its eigenvalues lie in  $\{0,2\}$ and
$\Delta u(q)\leq 2m$.
\end{lemma}
\begin{proof}
At a zero of the nonnegative smooth function $u$, $\nabla u(q)=0$. Twice differentiating $\Sigma_{k}(u_{k})^{2}=4u$ in normal coordinates gives
$$
2\sum_{k} u_{ik}u_{jk}+2\sum_{k} u_{k} u_{kij}=4u_{ij}.
$$
At $q$ the second term vanishes, so $\mathcal{A}^{2}=2\mathcal{A}$, and the eigenvalue equation $\mu^{2}=2\mu$ proves the assertion.
\end{proof}

\begin{lemma}\label{lemma 2.3}
Every complete connected totally geodesic isometric immersion $x:N^{m}\to\mathbb{R}^{N}$ is a global isometry onto an affine $m$-plane.
\end{lemma}
\begin{proof}
Fix $p\in N$ and put $P=dx_{p}(T_{p}N)$. For any $v\in T_{p}N$, completeness makes the geodesic $\gamma_v(t)=\exp_{p}(tv)$ available for all $t$. Since $B=0$,
$$\frac{d^{2}}{dt^{2}}x(\gamma_{v}(t))=B(\dot\gamma_{v},\dot\gamma_{v})=0,$$
and hence
$$x(\exp_{p}v)=x(p)+dx_{p}(v).$$
By Hopf-Rinow, $\exp_{p}$ is surjective. Therefore the image is precisely $x(p)+P$. If $x(q_{1})=x(q_{2})$, choose $v_{i}$ with $\exp_{p}v_{i}=q_{i}$; the displayed equation implies $dx_{p}(v_{1})=dx_{p}(v_{2})$, so $v_{1}=v_{2}$ and $q_{1}=q_{2}$. The bijective local isometry $x$ is a global isometry.
\end{proof}

\begin{lemma}\label{lemma 2.4}
If $(N,g)$ is complete and $f\in C^{1}(N)$ is bounded below, there are $p_{j}\in N$ such that
$$
f(p_{j})\to\inf_{N}f,\quad |\nabla f|(p_{j})\to 0.
$$
For a function bounded above, the analogous statement holds at the supremum.
\end{lemma}
\begin{proof}
Choose $\varepsilon_{j}\downarrow0$ and $q_{j}$ with $f(q_{j})\leq\inf f+\varepsilon^{2}_{j}$. Ekeland's variational principle \cite{Ekeland1974} on the complete metric space $(N,d)$ gives $p_{j}$ satisfying $f(p_{j})\leq f(q_{j})$ and
$$
f(q)\ge f(p_{j})-\varepsilon_{j} d(p_{j},q)\quad(\forall q\in N).
$$
Apply this to $q=\exp_{p_{j}}(tv)$ for unit $v$, first with $t>0$ and then with $t<0$, to find $|df_{p_{j}}(v)|\le\varepsilon_{j}$. Taking the supremum over $v$ gives the required gradient estimate. Replace $f$ by $-f$ for a maximum sequence.
\end{proof}

 \vskip10mm
\section{Rigidity of translators with constant $|\vec{H}|$}\label{section 3}

\begin{proposition}\label{proposition 3.1}
Let $x:M^{2}\to\mathbb{R}^{4}$ be a connected translator with nonzero fixed velocity $V$. If $|\vec{H}|^2$ is constant, then $\vec{H}=B=0$, even without completeness.
\end{proposition}
\begin{proof}
Set $T=V^{\top}$, $f=\langle V,x\rangle$ and $c=|\vec{H}|^2$. Because $V=T+\vec{H}$ is a fixed ambient vector, \eqref{2.2} yields, for all tangent $X$,
\begin{equation}\label{3.1}
\nabla_{X}T=A_{\vec{H}}X,\quad \nabla^{\perp}_{X}\vec{H}=-B(X,T).
\end{equation}
Taking the gradient and Laplacian of $f$ gives
\begin{equation}\label{3.2}
\nabla f=T,\quad |T|^2=|V|^2-c=a^2,\quad
\Delta f=\langle V,\vec{H}\rangle=c,
\end{equation}
where $a\geq0$. We will prove by contradiction that $c=0$. Assume $c>0$.

If $a=0$, then $f$ is constant, contrary to $\Delta f=c>0$. Thus $a>0$. The constant length of $T=\nabla f$ implies
$$
0=\frac{1}{2}\nabla|T|^2=\nabla^{2}f(T,\cdot)=\langle A_{\vec{H}}T,\cdot\rangle,
$$
so $A_{\vec{H}}T=0$. On a local oriented chart let $e_{1}=\frac{T}{a}$ and choose $e_{2}$ perpendicular to it. Since $\operatorname{tr}A_{\vec{H}}=\langle \vec{H},\vec{H}\rangle=c$, we obtain
\begin{equation}\label{3.3}
A_{\vec{H}}e_{1}=0,\quad A_{\vec{H}}e_{2}=ce_{2}.
\end{equation}
Put $k=\frac{c}{a}>0$. The first identity of \eqref{3.1} and metric compatibility give

\begin{equation}\label{3.4}
\nabla_{e_{1}}e_{1}=\nabla_{e_{1}}e_{2}=0, \quad \nabla_{e_{2}}e_{1}=ke_{2}, \quad \nabla_{e_{2}}e_{2}=-ke_{1}.
\end{equation}
Since $k$ is constant, these connection coefficients give
\begin{equation}\label{3.5}
[e_{1},e_{2}]=-ke_{2},\quad R(e_{1},e_{2})e_{2}=-k^{2}e_{1},\quad K=-k^{2}.
\end{equation}
The negative value is not yet a contradiction; we must impose the normal Codazzi equations.

Let $h=\sqrt c$ and choose a local orthonormal normal frame $\nu_{1}=\frac{\vec{H}}{h}$, $\nu_{2}$. Equation \eqref{3.3} yields
$\langle B_{11},\nu_{1}\rangle=\langle B_{12},\nu_{1}\rangle=0$ and
$\langle B_{22},\nu_{1}\rangle=h$. Since $B_{11}+B_{22}=h\nu_{1}$, for smooth functions $b,d$ we have
\begin{equation}\label{3.6}
B_{11}=b\nu_{2},\qquad B_{12}=d\nu_{2},\quad B_{22}=h\nu_{1}-b\nu_{2}.
\end{equation}
The Gauss equation gives
$$K=\langle b\nu_{2},h\nu_{1}-b\nu_{2}\rangle-|d\nu_{2}|^{2}=-b^{2}-d^{2}.$$
Together with \eqref{3.5}, this gives the fixed radius relation
\begin{equation}\label{3.7}
b^2+d^2=k^2.
\end{equation}
Define normal connection functions by
$\nabla^{\perp}_{e_{i}}\nu_{1}=\omega_{i}\nu_{2}$ and
$\nabla^{\perp}_{e_{i}}\nu_{2}=-\omega_{i}\nu_{1}$.
The second identity of \eqref{3.1} with $T=ae_{1}$ and \eqref{3.6} yield
\begin{equation}\label{3.8}
h\omega_{1}=-ab,\quad h\omega_{2}=-ad.
\end{equation}

\smallskip\noindent\emph{Step 3.}
Use \eqref{3.4} to expand the first Codazzi equation
$$
(\nabla_{e_1}B)(e_{2},e_{2})=(\nabla_{e_{2}}B)(e_{1},e_{2}).
$$
The $\nu_{2}$ component of the left-hand side is $h\omega_{1}-e_{1}(b)$. The right-hand side has $\nu_{2}$ component $e_{2}(d)+2kb$: indeed,
$$
(\nabla_{e_{2}}B)(e_{1},e_{2})
=\nabla^{\perp}_{e_{2}}B_{12}-kB_{22}+kB_{11}.
$$
Thus
\begin{equation}\label{3.9}
h\omega_{1}-e_{1}(b)=e_{2}(d)+2kb.
\end{equation}
Similarly, the $\nu_{2}$ component of
$(\nabla_{e_{1}}B)(e_{1},e_{2})=(\nabla_{e_{2}}B)(e_{1},e_{1})$ gives
\begin{equation}\label{3.10}
e_{1}(d)=e_{2}(b)-2kd.
\end{equation}
Substituting \eqref{3.8} gives the closed system
\begin{equation}\label{3.11}
\begin{cases}
e_{1}(b)+e_{2}(d)=-(a+2k)b,\\
e_{1}(d)-e_{2}(b)=-2kd.
\end{cases}
\end{equation}

Since $k>0$, equation \eqref{3.7} allows a smooth local angle $\theta$ with
$b=k\cos\theta$, $d=k\sin\theta$. Substituting into \eqref{3.11} and inverting the rotation matrix gives
\begin{equation}\label{3.12}
e_{1}(\theta)=a\sin\theta\cos\theta,\quad
e_{2}(\theta)=-2k-a\cos^{2}\theta.
\end{equation}
From \eqref{3.5}, $[e_{1},e_{2}](\theta)=-k e_{2}(\theta)$. Differentiating \eqref{3.12} gives
\begin{equation*}
\begin{aligned}
&e_{1}(e_{2}(\theta))=2a^{2}\sin^{2}\theta\cos^{2}\theta,\\
&e_{2}(e_{1}(\theta))=a(\cos^{2}\theta-\sin^{2}\theta)(-2k-a\cos^{2}\theta).
\end{aligned}
\end{equation*}
The commutator equality therefore reduces, after collecting terms, to
\begin{equation}\label{3.13}
a(a+3k)\cos^{2}\theta=2k(a+k).
\end{equation}
All coefficients are positive constants, so $\cos^{2}\theta$ must be constant. Locally, smoothness then forces $\theta$ to be constant, because each level set of $\cos^2$ consists of isolated points. But \eqref{3.12} gives
$e_{2}(\theta)=-2k-a\cos^{2}\theta<0$. This is a contradiction. Hence $c>0$ is impossible.

Now $\vec{H}=0$ and $V^{\perp}=0$, so $V=T$ is a nowhere-zero tangent vector, and \eqref{3.1} gives $B(X,V)=0$ for every $X$. Choosing $e_{1}=V/|V|$ yields $B_{11}=B_{12}=0$. Since $\vec{H}=B_{11}+B_{22}=0$, we also have $B_{22}=0$. Thus $B\equiv0$.
\end{proof}

\begin{proof}[Proof of Theorem \ref{theorem 1.1}]
Proposition \ref{proposition 3.1} gives $B=\vec{H}=0$, and Lemma \ref{lemma 2.3} identifies the complete immersion with an entire affine plane $p_{0}+P$.
Since $V$ is tangent to the plane at every point, it lies in $P$.
\end{proof}

 \vskip10mm
\section{Rigidity of self-expanders with constant $|\vec{H}|$}\label{section 4}
\begin{proof}[Proof of Theorem \ref{theorem 1.2}]
Write $c=|\vec{H}|^{2}\geq0$ and $T=x^{\top}$. For the constant $c$, the self-expander equation gives
\begin{equation}\label{4.1}
u=|x|^2-c=|T|^2\geq0.
\end{equation}
Directly from \eqref{2.1},
\begin{equation}\label{4.2}
\nabla u=2T,\quad |\nabla u|^2=4u,\quad
\Delta u=2\langle x,\vec{H}\rangle+4=4+2c.
\end{equation}
Lemma \ref{lemma 2.1} gives $q\in M$ with $u(q)=0$.
By Lemma \ref{lemma 2.2}, $\Delta u(q)\leq4$.
Equation~\eqref{4.2} therefore forces $c=0$ and $\vec{H}\equiv0$. Then equation \eqref{1.2} gives $x^{\perp}=0$, i.e. $x=T$ is tangent everywhere. Taking normal parts of $D_{X}x=X$ yields
\begin{equation}\label{4.3}
B(X,T)=0.
\end{equation}
At any point where $T\neq0$, choose $e_{1}=T/|T|$. Equation \eqref{4.3} shows $B_{11}=B_{12}=0$, and $\vec{H}=0$ gives $B_{22}=0$. By \eqref{4.2}, $\Delta u=4$, so $u$ cannot vanish on an open set and the set $\{T\neq0\}$ is dense. Continuity gives $B=0$ globally. By Lemma \ref{lemma 2.3}, $x$ parametrizes an entire affine plane. Finally $u(q)=0$ and $c=0$ imply $x(q)=0$, so the plane passes through the origin.
\end{proof}

 \vskip10mm
\section{Lagrangian structure and drift identities}\label{section 5}
For the remainder of the paper assume that $x:M^{2}\to\mathbb{C}^{2}$ is Lagrangian. Set
\begin{equation}\label{5.1}
C(X,Y,W)=\langle B(X,Y),JW\rangle
\end{equation}
Because $J$ is parallel and $J(TM)=NM$, the cubic tensor $C$ is totally symmetric. The covariant derivative $\nabla C$ is totally symmetric in its four arguments by the Codazzi equation. The normal curvature is identified with the tangent curvature through $J$. That is, $R^{\perp}(X,Y)JW=JR(X,Y)W$. Then $K=0$ implies flat normal connection.

The equations \eqref{1.1} and \eqref{1.2} can be written as

\begin{equation}\label{5.2}
\vec{H}=(\lambda x+V)^{\perp},
\quad (\lambda,V)=(0,V),\ V\neq0,\quad\text{or}\quad(1,0).
\end{equation}
Define
\begin{equation}\label{5.3}
Z=-(\lambda x+V)^{\top},\quad \rho=|\vec{H}|^2,\quad S=|B|^2,
\quad \mathcal{L}=\Delta-\langle Z,\nabla\cdot\rangle.
\end{equation}
Thus $\mathcal{L}=\Delta+\langle V^{\top},\nabla\cdot\rangle$ for translators and $\mathcal{L}=\Delta+\langle x^{\top},\nabla\cdot\rangle$ for self-expanders.  Direct differentiation of \eqref{5.2} gives
\begin{equation}\label{5.4}
\nabla_{X}^{\perp} \vec{H}=B(X,Z),\quad
\nabla_{X} Z=-\lambda X-A_{\vec{H}}X,\quad
\nabla\rho=2A_{\vec{H}}Z.
\end{equation}
The final formula follows from $$X(\rho)=2\langle\nabla_{X}^{\perp} \vec{H},\vec{H}\rangle=2\langle A_{\vec{H}}Z,X\rangle.$$
With the conventions in \eqref{5.2}-\eqref{5.3}, we obtain the following drift identities.

\begin{lemma}\label{lemma 5.1} For a Lagrangian surface satisfying \eqref{5.2}, we have
\begin{equation}\label{5.5}
\frac{1}{2}\mathcal{L}\rho=|\nabla^{\perp} \vec{H}|^{2}-\lambda\rho-|A_{\vec{H}}|^{2},
\end{equation}
and
\begin{equation}\label{5.6}
\frac{1}{2}\mathcal{L}S=|\nabla^{\perp}B|^{2}-\lambda S-\frac{3}{2}S^{2}
+2\rho S-\frac{1}{2}\rho^{2}-|A_{\vec{H}}|^{2}.
\end{equation}
Consequently, as $2K=\rho-S$,
\begin{equation}\label{5.7}
\mathcal{L}K=|\nabla^{\perp}\vec{H}|^{2}-|\nabla^{\perp}B|^{2}+(-2\lambda-2\rho+6K)K.
\end{equation}
\end{lemma}
\begin{proof}
Tracing the derivative of the first identity in \eqref{5.4} and using Codazzi gives
$$
\Delta^{\perp}\vec{H}=\nabla^{\perp}_Z\vec{H}-\lambda\vec{H}-\sum_{i}B(e_{i},A_{\vec{H}}e_{i}).
$$
Taking its inner product with $\vec{H}$ and using $\frac{1}{2}\Delta\rho=|\nabla^{\perp}\vec{H}|^{2}+\langle\Delta^{\perp}\vec{H},\vec{H}\rangle$ yields \eqref{5.5}.
For clarity, here is the algebraic part of the Simons formula used in \eqref{5.6}. With respect to a tangent orthonormal frame write
\begin{equation*}
\begin{aligned}
A_{Je_{1}}=\begin{bmatrix}
a & b \\
b & c
\end{bmatrix}, \quad
A_{Je_{2}}=\begin{bmatrix}
b & c \\
c & d
\end{bmatrix}.
\end{aligned}
\end{equation*}
Then
\begin{equation}\label{5.8}
S=a^{2}+3b^{2}+3c^{2}+d^{2},\quad
\rho=(a+c)^{2}+(b+d)^{2},\quad
A_{\vec{H}}=(a+c)A_{Je_{1}}+(b+d)A_{Je_{2}}.
\end{equation}
A direct expansion of these two symmetric matrices verifies
\begin{equation}\label{5.9}
\sum_{i,j=1}^{2}
\Big[(\operatorname{tr}(A_{Je_{i}} A_{Je_{j}}))^{2}+
|[A_{Je_{i}},A_{Je_{j}}]|^{2}\Big]=\frac{3}{2}S^{2}-2\rho S+\frac{1}{2}\rho^{2}+|A_{\vec{H}}|^{2},
\end{equation}
The standard Euclidean Simons identity for the second fundamental form, combined with \eqref{5.4}, is
$$
\frac{1}{2}\mathcal{L} S=|\nabla^{\perp}B|^{2}-\lambda S
-\sum_{i,j=1}^{2}
\Big[(\operatorname{tr}(A_{Je_{i}} A_{Je_{j}}))^{2}+
|[A_{Je_{i}},A_{Je_{j}}]|^{2}\Big].
$$
Substituting \eqref{5.9} into the Simons identity gives \eqref{5.6}. Finally, $2K=\rho-S$ and subtraction of the first two equations give \eqref{5.7}.
\end{proof}

\vskip10mm
\section{The cubic-tensor identity at negative constant curvature}\label{section 6}
The next calculation is the key point in excluding $K<0$. It is important that a frame identity valid only at one point is not differentiated as though it held on a neighbourhood.

\begin{lemma}\label{lemma 6.1}
Suppose $K=-\kappa<0$ is constant and \eqref{5.2} holds. At any point where $\rho>0$, choose an orthonormal tangent frame with
$\vec{H}=\sqrt{\rho}\,Je_{1}$ and write $Z=\xi e_{1}+\zeta e_{2}$. Then
\begin{equation}\label{6.1}
\kappa|Z|^{2}=\rho\xi^{2}+(\rho+8\kappa)(\lambda+\rho+3\kappa).
\end{equation}
\end{lemma}

\begin{proof}
The cubic tensor $C_{ijk}=\langle B(e_{i},e_{j}), Je_{k}\rangle$ and its covariant derivative are totally symmetric. At the point in question write
\begin{equation}\label{6.2}
(C_{111},C_{112},C_{122},C_{222})=(a,b,c,-b),\quad
h_{0}=a+c=\sqrt{\rho}.
\end{equation}
The last component $C_{222}$ is $-b$ because the $Je_{2}$ component of $\vec{H}$ is zero at this point. Directly from the Gauss equation,
\begin{equation}\label{6.3}
K=ac-c^{2}-2b^{2},\qquad S=a^{2}+3c^{2}+4b^{2}.
\end{equation}
For clarity, the three normal-valued coefficients are
$$
B_{11}=aJe_{1}+bJe_{2},\quad B_{12}=bJe_{1}+cJe_{2},
\quad B_{22}=cJe_{1}-bJe_{2}.
$$
At the point, denote the independent components of the symmetric four-tensor $\nabla C$ by
$$
p=C_{111,1},\quad q=C_{111,2},\quad
r=C_{112,2},\quad s=C_{122,2},\quad t=C_{222,2}.
$$
Because $\nabla^{\perp}_{X}\vec{H}=B(X,Z)$ by \eqref{5.4}, we have
\begin{equation}\label{6.4}
p+r=A=a\xi+b\zeta,\quad
q+s=F=b\xi+c\zeta,\quad
r+t=E=c\xi-b\zeta.
\end{equation}
Define
\begin{equation}\label{6.5}
P=r-\frac{h_{0}\xi}{8},\qquad Q=q-\frac{F}{2}.
\end{equation}
The multiplicities of symmetric entries yield
\begin{equation*}
|\nabla^{\perp}B|^{2}=p^{2}+4q^{2}+6r^{2}+4s^{2}+t^{2},\quad
|\nabla^{\perp}\vec{H}|^{2}=A^{2}+2F^{2}+E^{2}.
\end{equation*}
Using $p=A-r$, $s=F-q$, $t=E-r$,
\begin{equation}\label{6.6}
\begin{aligned}
|\nabla^{\perp}\vec{H}|^{2}-|\nabla^{\perp}B|^{2}
&=2r(A+E)-8r^{2}-8(q-F/2)^{2}\\
&=\frac{\rho\xi^{2}}{8}-8(P^{2}+Q^{2}),
\end{aligned}
\end{equation}
where $A+E=h_{0}\xi$.

One must use the differentiated full cubic invariant, rather than differentiating $C_{222}=-b$ on a neighborhood. In an arbitrary frame, the intrinsic Gauss expression is
$$
K=C_{111}C_{122}+C_{112}C_{222}-C_{112}^{2}-C_{122}^{2}.
$$
Differentiate this invariant formula first, and only then evaluate in the adapted frame where $C_{222}=-b$. Total symmetry of $\nabla C$ gives
\begin{equation}\label{6.7}
K_{1}=cp+(a-2c)r-3bq+bs,
\end{equation}
\begin{equation}\label{6.8}
K_{2}=cq+(a-2c)s-3br+bt.
\end{equation}
Inserting \eqref{6.4} and \eqref{6.5} simplifies these to
\begin{equation}\label{6.9}
K_{1}=\Big(\frac{\rho}{8}+\frac{K}{2}\Big)\xi+(h_{0}-4c)P-4bQ,
\end{equation}
\begin{equation}\label{6.10}
K_{2}=\frac{K}{2}\zeta-4bP-(h_{0}-4c)Q.
\end{equation}
Let
\begin{equation*}
\mathcal{B}=\begin{bmatrix}
h_{0}-4c & -4b \\
-4b & -h_{0}+4c
\end{bmatrix},\quad D=\rho-8K.
\end{equation*}
From $h_{0}=a+c$ and \eqref{6.3},
$$
(h_{0}-4c)^{2}+16b^{2}=(a-3c)^{2}+16b^{2}
=(a+c)^{2}-8(ac-c^{2}-2b^{2})=D,
$$
so $\mathcal{B}^{2}=DI$. Since $K$ is constant and $D=\rho+8\kappa>0$, equations \eqref{6.9}-\eqref{6.10} imply
\begin{equation*}
\mathcal{B}\begin{pmatrix}
P\\
Q
\end{pmatrix}
=
-\begin{pmatrix}
(\rho/8+K/2)\xi\\
(K/2)\zeta
\end{pmatrix}
\end{equation*}
Taking norms and using $\mathcal{B}^{2}=DI$ yields
\begin{equation}\label{6.11}
P^{2}+Q^{2}=\frac{1}{D}\Big[
  \big(\frac{\rho}{8}+\frac{K}{2}\big)^{2}\xi^{2}+\frac{K^{2}}{4}\zeta^{2}\Big].
\end{equation}
Combining \eqref{6.6} and \eqref{6.11}, and simplifying the coefficient of $\xi^{2}$, gives the full identity
\begin{equation}\label{6.12}
|\nabla^{\perp}\vec{H}|^{2}-|\nabla^{\perp}B|^{2}
=-\frac{2K^{2}}{D}|Z|^{2}-\frac{2K\rho}{D}\xi^{2}.
\end{equation}
Indeed,
$$\frac{\rho(\rho-8K)-(\rho+4K)^{2}}{8D}=-\frac{2K(\rho+K)}{D}.$$
Finally, since $\mathcal{L}K=0$, inserting \eqref{6.12} into \eqref{5.7}, setting $K=-\kappa$, and multiplying by $D/(2\kappa)>0$ gives
$$
0=-\kappa|Z|^{2}+\rho\xi^{2}+(\rho+8\kappa)(\lambda+\rho+3\kappa),
$$
which is \eqref{6.1}.
\end{proof}

\vskip10mm
\section{Exclusion of nonzero constant Gaussian curvature}\label{section 7}
\begin{proposition}\label{proposition 7.1}
No complete immersed translator or self-expander in $\mathbb{R}^{4}$ has positive constant Gaussian curvature.
\end{proposition}
\begin{proof}
For a surface, $\operatorname{Ric}=Kg$. A positive constant $K$ makes $M$ compact by the Bonnet-Myers theorem.
For a translator, $f=\langle V,x\rangle$ satisfies $\Delta f=\langle V,\vec{H}\rangle=|\vec{H}|^{2}$. Integration gives $\vec{H}=0$, then $\Delta|x|^{2}=4$, contradicting compactness. For a  self-expander,
$$
\frac{1}{2}\Delta(|x|^{2})=2+\langle x,\vec{H}\rangle=2+|\vec{H}|^{2}>0,
$$
which is also impossible on a compact manifold.
\end{proof}

\begin{proposition}\label{proposition 7.2}
A complete Lagrangian self-expander in $\mathbb{C}^{2}$ cannot have negative constant Gaussian curvature.
\end{proposition}
\begin{proof}
Suppose $K=-\kappa<0$. Let $f=|x|^{2}/2\geq0$. For self-expanders, $Z=-x^{\top}=-\nabla f$ and $\lambda=1$. On $\Omega=\{\rho>0\}$, Lemma \ref{lemma 6.1} implies
\begin{equation}\label{7.1}
\kappa|\nabla f|^2=\rho\xi^{2}+(\rho+8\kappa)(1+\rho+3\kappa)
\geq8\kappa(1+3\kappa).
\end{equation}
We show that $\Omega$ is dense. If $\vec{H}=0$ on a nonempty open set, the self-expander equation gives $x^{\perp}=0$ there. For tangent vectors $X,Y$, the second equation of \eqref{5.4} yields
$$
\nabla^{2}f(X,Y)=\langle X,Y\rangle+\langle x^{\perp},B(X,Y)\rangle=\langle X,Y\rangle.
$$
Thus $\nabla_{X}\nabla f=X$ there. Commuting these derivatives yields $R(X,Y)\nabla f=0$. For a $2$-dimensional Riemannian manifold, its curvature tensor is completely determined by the Gaussian curvature:
$$
R(X,Y)W = K(\langle W,Y\rangle X-\langle W,X\rangle Y)
$$
where the corresponding curvature sign convention is adopted. Choose a local orthonormal frame $\{e_{1}, e_{2}\}$, and write
$\nabla f=f_{1}e_{1} + f_2e_{2}$. It follows from $R(X,Y)\nabla f=0$ and $K=-\kappa \neq 0$ that
$f_{1}=f_{2}=0$. Then on a surface of nonzero Gaussian curvature this implies $\nabla f=0$, contrary to $\nabla^{2}f=g$.
Consequently the continuous estimate \eqref{7.1} holds on all of $M$.
Since $f$ is bounded below, Lemma \ref{lemma 2.4} supplies points $p_{j}$ with $|\nabla f|(p_{j})\to0$. Inequality \eqref{7.1} gives instead $|\nabla f|^{2}\geq8(1+3\kappa)>0$, a contradiction.
\end{proof}

\begin{proposition}\label{proposition 7.3}
A complete Lagrangian translator in $\mathbb{C}^2$ cannot have negative constant Gaussian curvature.
\end{proposition}
\begin{proof}
Suppose $K=-\kappa<0$. 
The translator equation and Gauss identity yield
\begin{equation}\label{7.2}
0\leq\rho=|V^{\perp}|^{2}\leq|V|^{2},\quad
|Z|=|V^{\top}|\leq|V|,\quad
S=\rho+2\kappa\leq|V|^{2}+2\kappa.
\end{equation}
If $\vec{H}=0$ identically, then $V$ is tangent, and from $D_{X}V=0$, it is a nonzero parallel tangent field which forces $K=0$, contrary to hypothesis. Therefore $\rho_{\ast}=\sup_{M}\rho>0$.
Apply Lemma \ref{lemma 2.4} to $-\rho$. We find $p_{j}$ with
\begin{equation}\label{7.3}
\rho(p_{j})\to\rho_{\ast},\quad |\nabla\rho|(p_{j})\to0.
\end{equation}
For large $j$, use the mean curvature frame of Lemma \ref{lemma 6.1}, so that $\vec{H}=\sqrt\rho Je_{1}$.
By \eqref{7.2}, the coefficients $a,b,c$ of \eqref{6.2} and $\xi,\zeta$ of $Z$ are uniformly bounded.
After taking a subsequence, these five real numbers converge to
$a,b,c,\xi,\zeta$, with
$h_{0}=a+c=\sqrt{\rho_{\ast}}>0$. Let
\begin{equation*}
\mathcal{A}=\begin{bmatrix}
a & b \\
b & c
\end{bmatrix},\quad
z=\begin{pmatrix}
\xi\\
\zeta
\end{pmatrix}
\end{equation*}
At the points $p_{j}$, the gradient identity in \eqref{5.4} is $\nabla\rho=2A_{\vec{H}}Z$. In the chosen frame $A_{\vec{H}}=h_{0}\mathcal{A}$. Passing to the limit in \eqref{7.3} gives
\begin{equation}\label{7.4}
\mathcal{A}z=0.
\end{equation}
Similarly, Lemma \ref{lemma 6.1} with $\lambda=0$ passes to the limit as
\begin{equation}\label{7.5}
\kappa|z|^{2}=\rho_{\ast}\xi^{2}+(\rho_{\ast}+8\kappa)(\rho_{\ast}+3\kappa).
\end{equation}

The right-hand side is positive, so $z\neq0$. Equation \eqref{7.4} therefore implies that $\det \mathcal{A}=ac-b^{2}=0$. Combining this with the Gauss relation \eqref{6.3} gives
\begin{equation}\label{7.6}
b^{2}+c^{2}=\kappa,\quad h_{0}c=(a+c)c=\kappa.
\end{equation}
The eigenvalues of the symmetric $2\times2$ matrix $\mathcal{A}$ are $0$ and $h_{0}$. Hence $\ker \mathcal{A}$ is one-dimensional.
Direct calculation yields $\mathcal{A}\upsilon=0$ for $\upsilon=(c,-b)^{\top}$. Then $\upsilon\in\ker \mathcal{A}$.
Since $z\neq0$ and $\mathcal{A}z=0$, there necessarily exists a non-zero real number $\tau$ such that
$z=\tau\upsilon$. Thus, $$\xi=\tau c, \quad \zeta=-\tau b, \quad |z|^{2}=\tau^{2}(c^{2}+b^{2})=\tau^{2} \kappa$$ from \eqref{7.6}.
Therefore, inserting $h_{0}c=\kappa$ and $\rho_{\ast}=h^{2}_{0}$ yields
\begin{equation}\label{7.7}
\rho_{\ast}\xi^{2}= h^{2}_{0} \tau^{2} c^{2}=\tau^{2} (h_{0}c)^{2}=\tau^{2} \kappa^{2}=\kappa|z|^{2}.
\end{equation}
Substitution of \eqref{7.7} into \eqref{7.5} forces
$0=(\rho_{\ast}+8\kappa)(\rho_{\ast}+3\kappa)>0$, the desired contradiction.
\end{proof}

From the above analysis, the following result can be derived.
\begin{corollary}\label{corollary 7.1}
Every complete connected Lagrangian translator or self-expander in $\mathbb{C}^{2}$ with constant Gaussian curvature is flat.
\end{corollary}

\vskip10mm
\section{Global splitting in the flat case}\label{section 8}
\begin{proposition}\label{proposition 8.1}
Let $x:M^{2}\to\mathbb{C}^{2}$ be a complete, connected, flat Lagrangian immersion satisfying \eqref{5.2}. Up to a unitary transformation, its lift to the universal cover is a Cartesian product of two complete, unit-speed planar soliton curves in orthogonal complex lines.
\end{proposition}
\begin{proof}
As $J$ identifies the tangent and normal connections, $K=0$ implies the normal curvature $R^{\bot}$ vanishes. By the Ricci equation all shape operators commute. At a point with $B\neq0$, at least one shape operator is nonscalar: otherwise Lagrangian symmetry would force $B=0$. A nonscalar self-adjoint operator in dimension two has distinct eigenvalues, hence determines a smooth local orthonormal eigenframe $\{e_{1},e_{2}\}$; commutation diagonalizes the remaining shape operators simultaneously. In this common principal frame, for any normal vector
$\nu \in N_{p}M$, we have $\langle B(e_{1}, e_{2}),\nu\rangle=0$. That is, $B(e_{1}, e_{2})=0$.
Total symmetry of $C$ then yields
\begin{equation}\label{8.1}
B(e_{1},e_{1})=k_{1}Je_{1},\quad
B(e_{2},e_{2})=k_{2}Je_{2}.
\end{equation}
Write $\nabla_{e_i}e_{1}=\omega_{i} e_{2}$, $\nabla_{e_{i}}e_{2}=-\omega_{i}e_{1}$. In the diagonal frame, total symmetry of $\nabla C$ gives
\begin{equation}\label{8.2}
e_{2}(k_{1})=\omega_{1}k_{1},\quad e_{1}(k_{2})=-\omega_{2}k_{2},
\quad k_{1}\omega_{2}+k_{2}\omega_{1}=0.
\end{equation}
For completeness, these are obtained by equating $C_{111,2}$ with $C_{112,1}$, $C_{222,1}$ with $C_{122,2}$, and $C_{112,2}$ with $C_{122,1}$, respectively.
Moreover, the $Je_{2}$ component of $\nabla^{\bot}_{e_{1}}\vec{H}=B(e_{1},Z)$ equals zero. Since $\vec{H}=k_{1}Je_{1}+k_2Je_{2}$, this component is $e_{1}(k_{2})+k_{1}\omega_{1}$. Together with \eqref{8.2},
$$
k_{1}\omega_{1}-k_{2}\omega_{2}=0,\quad
k_{2}\omega_{1}+k_{1}\omega_{2}=0.
$$
The determinant is $k^{2}_{1}+k^{2}_{2}>0$ at every point of this neighborhood. Thus
\begin{equation}\label{8.3}
\omega_{1}=\omega_{2}=0
\end{equation}
on this nonflat neighbourhood. In particular, the local diagonalising frame is parallel.

The universal Riemannian cover $(\tilde{M}, \tilde{g})$ is complete, simply connected, and flat; hence it is globally isometric to the Euclidean plane with coordinates $(s,t)$. Extend the local parallel eigenframe uniquely to global parallel fields $\partial_{s},\partial_{t}$ on the cover. As $\nabla_{\partial_{s}}\partial_{t}=0$ and $B(\partial_{s},\partial_{t})=0$ on the initial neighbourhood,
\begin{equation}\label{8.4}
x_{st}=\nabla_{\partial_{s}}\partial_{t}+B(\partial_{s},\partial_{t})=0
\end{equation}
To propagate this identity, use the soliton equation in the global Euclidean coordinates:
\begin{equation}\label{8.5}
x_{ss}+x_{tt}=\lambda x+V
 -\langle\lambda x+V,x_{s}\rangle x_{s}
 -\langle\lambda x+V,x_{t}\rangle x_{t}.
\end{equation}
Here $|x_{s}|=|x_{t}|=1$ and $\langle x_{s},x_{t}\rangle=0$. Equation \eqref{8.5} is a semilinear elliptic system with analytic nonlinearities; standard analytic elliptic regularity \cite{Morrey2008} implies that $x$ is real analytic. Therefore the analytic tensor $x_{st}$, which vanishes on the open set in \eqref{8.4}, vanishes throughout $\mathbb{R}^2$. If $B\equiv0$ instead, the same conclusion follows immediately. Thus
\begin{equation}\label{8.6}
x_{st}=0 \quad \text{on }\mathbb{R}^{2}.
\end{equation}
Integrating gives $x(s,t)=\alpha(s)+\beta(t)$,
with $|\alpha'|=|\beta'|=1$. The orthogonality of $x_{s}, x_{t}$ and the Lagrangian condition yield
\begin{equation}\label{8.7}
\langle\alpha^{'},\beta^{'}\rangle=\langle J\alpha^{'},\beta^{'}\rangle=0.
\end{equation}
By \eqref{8.1} on the initial neighbourhood and analyticity,
$$
\alpha^{''}(s)\in\operatorname{span}\{J\alpha^{'}(s)\},\quad
\beta^{''}(t)\in\operatorname{span}\{J\beta^{'}(t)\}.
$$
Therefore the real two-planes $L_{1}=\operatorname{span}\{\alpha^{'},J\alpha^{'}\}$ and $L_{2}=\operatorname{span}\{\beta^{'},J\beta^{'}\}$ are each independent of their parameters. Equation \eqref{8.7} gives $L_{1}\perp L_{2}$. They are fixed orthogonal complex lines in $\mathbb{C}^2$. Choosing a unitary $U$ that sends the coordinate lines to $L_{1},L_{2}$ and absorbing the constant components of $x$ into the factors, we obtain
$$
x(s,t)=U(\gamma_{1}(s),\gamma_{2}(t)).
$$
The normal projection in \eqref{5.2} splits over these complex lines, giving
$$
\gamma^{''}_{i}=(V_i)^{\perp_{\gamma_{i}}}\ \text{for translators},\quad
\gamma^{''}_{i}=(\gamma_{i})^{\perp_{\gamma_{i}}}\ \text{for expanders},
$$
where $U^{-1}(\lambda x+V)=(\lambda\gamma_{1}+V_{1}, \lambda\gamma_{2}+V_{2})$ and $\perp_{\gamma_{i}}$ denotes the projection onto the normal direction of the corresponding plane curve.
Both parameters traverse all of $\mathbb{R}$, so the curve factors are complete.
\end{proof}

\vskip10mm
\section{Planar curve factors and proofs of classification}\label{section 9}
\subsection{Translating plane curves}
Let $\gamma:\mathbb{R}\to\mathbb{C}$ be complete and unit speed with $\gamma^{''}=v^{\perp_{\gamma}}$ for fixed $v\in\mathbb{C}$. If $v=0$, the curve $\gamma$ is a line. For $v\neq0$, rotate the plane so that $v=(0,a)$ with $a=|v|>0$. Write $T=\gamma^{'}=(\cos\theta,\sin\theta)$ and $N=JT=(-\sin\theta,\cos\theta)$. The curvature equation is
\begin{equation}\label{9.1}
\theta^{'}=\langle v,N\rangle=a\cos\theta.
\end{equation}
An equilibrium solution is a straight line parallel to $v$. Every nonequilibrium solution has $\cos\theta$ of a fixed sign; reversing orientation allows $-\pi/2<\theta<\pi/2$. On that branch the curve is a graph $y=y(u)$, and \eqref{9.1} becomes
\begin{equation}\label{9.2}
\frac{y^{''}}{1+(y^{'})^{2}}=a.
\end{equation}
Integrating twice yields, after translations,
\begin{equation}\label{9.3}
\gamma_{a}(u)=\Big(u,-\frac{1}{a}\log\cos(au)\Big),\qquad
-\frac{\pi}{2a}<u<\frac{\pi}{2a}.
\end{equation}

The arclength form $ds=\sec(au)du$ diverges at both endpoints, so the maximal arclength parametrisation has domain $\mathbb{R}$. It is injective because the horizontal coordinate is strictly monotone, and proper because its height diverges at both ends. Uniqueness of the curve ODE shows that a complete curve must be one of these maximal solutions. For a prescribed velocity $V$, each factor velocity $V_{i}=(U^{-1}V)_{i}$ is fixed: a line must be parallel to $V_{i}$ when $V_{i}\neq0$, and a grim reaper has speed parameter $|V_{i}|$.

\subsection{Expanding plane curves}
Let $\gamma:\mathbb{R}\to\mathbb{C}$ be a unit-speed solution of $\gamma^{''}=\gamma^{\perp_{\gamma}}$. Put
$$
T=\gamma^{'},\quad N=JT,\quad p=\langle\gamma,T\rangle,
\quad k=\langle\gamma,N\rangle.
$$
Since $T^{'}=kN$ and $N^{'}=-kT$, differentiating gives
\begin{equation}\label{9.4}
p^{'}=1+k^{2},\quad k^{'}=-pk,\quad T^{'}=kN.
\end{equation}
As $p^{'}\geq1$, a solution defined for all $s\in\mathbb{R}$ has $p(s)\to+\infty$ as $s\to+\infty$ and $p(s)\to-\infty$ as $s\to-\infty$, so it has a unique zero. Translate the arc-length parameter to make this zero occur at $s=0$, reverse orientation if necessary, and rotate the ambient plane so that
\begin{equation}\label{9.5}
p(0)=0,\quad k(0)=a\geq0,
\quad T(0)=(1,0),\quad\gamma(0)=(0,a).
\end{equation}
If $a=0$, uniqueness gives $k=0$, $p=s$ and the line through the origin. Suppose $a>0$. By uniqueness $k$ remains positive, and on $s\geq0$,
$$
p(s)\geq s,\quad 0<k(s)\leq a,
\quad 1\leq p^{'}(s)\leq 1+a^{2}.
$$
The analogous backward bounds hold for $s\leq0$. Therefore the solution of \eqref{9.4} exists for all $s\in\mathbb{R}$. Conversely, solving \eqref{9.4} and then putting $\gamma=pT+kN$ gives
$$
\gamma^{'}=(p^{'}-k^{2})T+(pk+k^{'})N=T,
$$
so the reconstructed curve satisfies the desired expander equation.

To prove embeddedness, write $T=(\cos\theta,\sin\theta)$ with $\theta(0)=0$; thus $\theta^{'}=k>0$. For $s>0$ set $\phi=\arctan(p/k)$. Direct differentiation yields
\begin{equation}\label{9.6}
\phi^{'}=\frac{k p^{'}-p k^{'}}{p^{2}+k^{2}}
   =k+\frac{k}{p^{2}+k^{2}}>k=\theta^{'}.
\end{equation}
Here $\phi(0)=\theta(0)=0$, hence
$0<\theta(s)<\phi(s)<\pi/2$ on $s>0$. The symmetry of the initial-value system implies $p$ and $\theta$ are odd and $k$ is even. Therefore
$-\pi/2<\theta(s)<\pi/2$ for all $s\in\mathbb{R}$. It follows that $\frac{d}{ds}\gamma_{1}(s)=\cos\theta(s)>0$, proving injectivity. Also $p(s)\geq s$ for $s\geq0$ and $p(s)\leq s$ for $s\leq 0$; since
$$
|\gamma(s)|^{2}=p(s)^{2}+k(s)^{2}\geq s^{2},
$$
the curve is proper. Denote the normalized curve by $\Gamma_{a}$, with $\Gamma_{0}=\mathbb{R}$. Every complete planar expanding curve has this form up to rotation and reparametrization. The planar classification was established in \cite{Halldorsson2012}.

\subsection{Descent to the surface and proof of the main theorems}
\begin{proof}[Proofs of Theorems \ref{theorem 1.3} and \ref{theorem 1.4}]
Corollary \ref{corollary 7.1} gives $K=0$. By Proposition \ref{proposition 8.1}, the universal cover $\tilde{M}\cong\mathbb{R}^2$ has lifted immersion
$\tilde{x}(s,t)=U(\gamma_{1}(s),\gamma_{2}(t))$.
The preceding subsections show that each complete factor is injective and proper, so the Cartesian product is injective and proper. If $\sigma$ is a deck transformation, then $\tilde{x}\circ\sigma=\tilde{x}$. Injectivity forces $\sigma$ to be the identity. Hence the deck group is trivial and $M$ itself is isometric to $\mathbb{R}^{2}$. The original immersion is the proper embedded product stated in the theorems.

Conversely, for complete unit-speed curves in orthogonal complex lines, the product has metric $ds^{2}+dt^{2}$ and is flat and complete. Its tangent planes are Lagrangian, and both its mean curvature vector and the normal projection of $V$ or $x$ split into the respective planar components. The factor curve equations imply the corresponding surface soliton equation.
\end{proof}

\noindent{\bf Acknowledgements.}
The first author was partially supported by NSFC grant No.12401060, Natural Science Foundation of Henan Province grant No.262300421232. The second author was partially supported by NSFC grant No. 12501063, Postdoctoral Fellowship Program of CPSF
grant No. GZC20252036.

\noindent{\bf Declarations}

\noindent{\bf Conflict of interest} There are no conflicts of interest with third parties.

\noindent{\bf Data availability} No datasets were generated or analysed in this study.

\end{document}